\documentclass{article}
\usepackage[english]{babel}
\usepackage{latexsym,amssymb}
\usepackage[cp1251]{inputenc}
\usepackage{amsfonts,amssymb}
\usepackage{euscript}
\newcounter{theorem} 
\newcounter{lemma} 

\begin{document}

\vspace*{7mm}

\Large

\thispagestyle{empty}

\Large
 \begin{center}
    \textbf{Vitalii Shpakivskyi\footnote{Institute of Mathematics of the NAS of Ukraine, Kyiv, Ukraine; shpakivskyi86@gmail.com}}
 \end{center}
\vspace{7mm}
\Huge
 \begin{center}
   \textbf{
PART 3}
 \end{center}
\vspace{3mm}

\LARGE
 \begin{center}
   \textbf{Algebras with Variable Operator-Valued Structural Constants and Their Applications}
 \end{center}
\vspace{5mm}






 
\Large

In Part 1, we considered monogenic functions taking values in operator algebras with ``constant'' structural operators. In Part 2, algebras with variable functional structural coefficients were considered. In this part, we generalize the construction introduced in Part 1 by allowing the structural operators to depend on the coordinates of the hypercomplex variable as well as on additional parameters.

\section{Basic Concepts}

In this paper, we use all the concepts and notation introduced in Part 1. We briefly recall them below.

Let $\mathbb{K}\in\{\mathbb{R},\mathbb{C}\}$.
Let $X$ be a linear topological space over the field $\mathbb{K}$. Examples of spaces $X$ were given in Part 1.

Denote by $\mathcal{A}(X)$ a commutative associative algebra of linear operators
mapping the space $X$ into itself. Examples of algebras $\mathcal{A}(X)$ were also given in Part 1.

Let $\mathfrak{A}_n(\mathcal{A}(X),\star)$ be an $n$-dimensional algebra with operator-valued structural constants.
The algebra $\mathfrak{A}_n(\mathcal{A}(X),\star)$ is not necessarily commutative or associative. The $\star$-multiplication in this algebra is defined by
\begin{equation}\label{tabl-mnozh-zagalna-zmin-oper-stare}
E_i\star E_j=\sum_{k=1}^{n}\Upsilon_{ij}^{k}E_k,
\qquad i,j=1,2,\ldots,n,
\end{equation}
where $\Upsilon_{ij}^{k}\in\mathcal{A}(X)$, $i,j,k=1,\ldots,n$, are certain operators called the operator-valued structural constants of the algebra.
Clearly, the operators $\Upsilon_{ij}^{k}$ may have ``variable'' coefficients (although the use of the term ``variable coefficients'' for an operator is not entirely appropriate).
Examples of such operators were given in Part 1: Examples 1, 2, 5, 6, and 12.2--12.4.
However, when we consider functions taking values in the algebra $\mathfrak{A}_n(\mathcal{A}(X),\star)$, the notion of variable coefficients acquires a precise meaning.

Let $D$ be a domain in $\mathbb{R}^m$, where $1\leq m\leq n$.
In what follows, we consider algebras $\mathfrak{A}_n(\mathcal{A}(X),\star)$ whose structural operators depend on a point
$\xi=(x_1,x_2,\ldots,x_m)$ of the domain $D$.

\textbf{Definition 1.}\label{def-operatorna-algebra-zmin-operator}
Let $X$ be a linear topological space over the field $\mathbb{K}$, and let $\mathcal{A}(X)$ be a commutative associative algebra of linear operators mapping $X$ into itself.
To each point $\xi=(x_1,x_2,\ldots,x_m)$ of the domain $D$, we associate an algebra
$\mathfrak{A}_n(\mathcal{A}(X),\star;\xi)$ according to the multiplication rule
(\ref{tabl-mnozh-zagalna-zmin-oper-stare}), in which the operators $\Upsilon_{ij}^{k}$ depend on
$x_1,x_2,\ldots,x_m$ and, possibly, on additional variables $t_1,t_2,\ldots,t_r$, that is,
$$
\Upsilon_{ij}^{k}
=
\Upsilon_{ij}^{k}(x_1,x_2,\ldots,x_m;t_1,t_2,\ldots,t_r)
\in\mathcal{A}(X).
$$
Thus, for a fixed collection of structural operators
$$
\Upsilon:=
\{\Upsilon_{ij}^{k}:D\rightarrow\mathcal{A}(X),\,
i,j,k=1,\ldots,n\},
$$
each point $\xi$ of the domain $D$ is associated with an algebra
$\mathfrak{A}_n(\mathcal{A}(X),\star;\xi)$.
The collection of all such algebras
$$
\mathfrak{A}_n(D,\Upsilon):=
\Biggr\{
\mathfrak{A}_n(\mathcal{A}(X),\star;\xi)
\Biggr\}_{\xi\in D}
$$
will be called a \textbf{family of operator algebras with variable structural operators} $\Upsilon$, or, briefly, an
\textbf{algebra with variable operator-valued structural constants}.

Examples of such families of algebras are provided by Examples 5, 6, and 12.2--12.4 in Part 1.

\vskip2mm
\textbf{Remark 1.}
The basis elements $E_1,E_2,\ldots,E_n$ of the algebra
$\mathfrak{A}_n(D,\Upsilon)$ do not depend on the variables
$x_1,x_2,\ldots,x_m,t_1,t_2,\ldots,t_r$.

The concept of a family of algebras with variable operator coefficients generalizes the classical concept of an algebra given in \cite{Pierce}.

\section{Monogenic Functions in the Commutative Algebra $\mathfrak{A}_2(\mathcal{A}(X),\star)$ with Variable Operator-Valued Structural Constants}

Next, we generalize elements of the theory of monogenic functions developed in the monograph \cite{Plaksa-Shpakivskyi-2023}
 to the case of a commutative associative algebra with variable structural constants.

We demonstrate our approach using a special two-dimensional commutative associative algebra with variable operator-valued structural constants.

As before, consider the two-dimensional commutative associative operator algebra $\mathfrak{A}_2(\mathcal{A}(X),\star)$ with basis
$E_1,E_2$ and multiplication table
\begin{equation}\label{tabl-mnozh-alg-2-zmin-operat}
\begin{tabular}{c||c|c|}
$\star$ & $E_1$ & $E_2$ \\
\hline \hline
$E_1$ & $E_1$ & $E_2$ \\
\hline
$E_2$ & $E_2$ & $AE_1+BE_2$ \\
\hline
\end{tabular}\,\,.
\end{equation}

Let now $D$ be a domain in $\mathbb{R}^2$ and $(x,y)\in D$.
We assume that the operators $A$ and $B$ in the multiplication table (\ref{tabl-mnozh-alg-2-zmin-operat})
depend on $x,y$ and, possibly, on additional variables $t_1,t_2,\ldots,t_r$, that is,
$$
A=A(x,y;t_1,t_2,\ldots,t_r) \in \mathcal{A}(X),\qquad
B=B(x,y;t_1,t_2,\ldots,t_r)\in \mathcal{A}(X).
$$
Thus, for a fixed set of structural operators
$\Upsilon=\{A,B:D\rightarrow\mathcal{A}(X)\}$, each point $(x,y)$ of the domain $D$
is associated with an algebra $\mathfrak{A}_2(\mathcal{A}(X),\star;x,y)$. In what follows, we denote the collection of all these algebras by
$$
\mathfrak{A}_2(D,\Upsilon):=
\Biggr\{\mathfrak{A}_2(\mathcal{A}(X),\star;x,y)\Biggr\}_{(x,y)\in D}.
$$
 
\vskip2mm
\textbf{Remark 2.} The basis elements $E_1,E_2$ are the same for the entire family of algebras $\mathfrak{A}_2(D,\Upsilon)$ and do not depend on the variables
$x,y;t_1,t_2,\ldots,t_r$.

All algebraic properties established in Part 1 (Theorems 1--8) remain valid for the family of algebras $\mathfrak{A}_2(D,\Upsilon)$.

\textbf{Definition 2.} For simplicity of exposition, the element
$Z=xE_1+yE_2,\,\, (x,y)\in D$,
will be called the \textbf{variable of the family of algebras} $\mathfrak{A}_2(D,\Upsilon)$. The element $dZ=dxE_1+dyE_2$ is called the differential of the variable $Z$.
Applications may, however, require consideration of other variables, for example,
$Z=xi_1+yi_2$, where $i_1,i_2$ are certain vectors of the algebra $\mathfrak{A}_2(D,\Upsilon)$ satisfying prescribed relations.

Consider a function $\Phi:D\to \mathfrak{A}_2(D,\Upsilon)$ which assigns to each point $(x,y)\in D$ an element of the algebra
$\Biggr\{\mathfrak{A}_2(\mathcal{A}(X),\star;x,y)\Biggr\}_{(x,y)\in D}$, that is,
$\Phi(x,y)\in \Biggr\{\mathfrak{A}_2(\mathcal{A}(X),\star;x,y)\Biggr\}_{(x,y)\in D}$.
Since all algebras in the family $\mathfrak{A}_2(D,\Upsilon)$ have the common basis $\{E_1,E_2\}$, the function $\Phi$ can be represented in the form
\begin{equation}\label{rozklad za basysom -2-vym-algebra-zi-zminnymy-operatoramy}
\Phi(Z)=U(x,y)E_1+V(x,y)E_2,
\end{equation}
where $U,V:D\to\mathcal{A}(X)$.
The functions $U,V$ will be called the components of the function $\Phi$.

\vskip2mm
\textbf{Definition 3.} As before, let $D$ be a domain in $\mathbb{R}^2$.
Let $U(x,y):D\to\mathcal{A}(X)$ be an operator-valued function depending on the variables $x,y$.
The partial derivatives of the operator-valued function $U$ are understood in the usual sense as derivatives with respect to the variables $(x,y)\in D$.
The partial derivatives of the $\mathfrak{A}_2$-valued function
$\Phi(Z)=U(x,y)E_1+V(x,y)E_2$ are understood as the partial derivatives of the corresponding $\mathcal{A}(X)$-valued components $U$ and $V$.

\vskip2mm
\textbf{Example 1.}
Let $X=C^\infty(\mathbb R)$ and $A=\frac{d}{dt}$.
Consider the operator-valued function $U(x,y)=x^2A+yI$.
Then
$$
U_x(x,y)=2xA, \qquad U_y(x,y)=I.
$$

\vskip2mm
\textbf{Example 2.}
Let $L=\frac{d}{dt}$ and $U(x,y)=xL+y^2L^2$.
Then
$$
U_x=L, \qquad U_y=2yL^2.
$$

\vskip2mm
\textbf{Example 3.}
Let the Volterra operator $\mathcal{V}$ be defined by
$$
(\mathcal Vf)(t)=\int_0^t f(\tau)\,d\tau.
$$
Set
$$
U(x,y)=e^x\mathcal V+\sin y\,\mathcal V^2.
$$
Then
$$
U_x=e^x\mathcal V,
\qquad
U_y=\cos y\,\mathcal V^2.
$$

\vskip2mm
\textbf{Example 4.}
Let the operator $A(x,y)$ act according to the rule
$$
(A(x,y)f)(t)=\bigl(x+ty\bigr)f'(t).
$$
Set $U(x,y)=A(x,y)$.
Then $(U_xf)(t)=f'(t)$, that is, $U_x=\frac{d}{dt}$,
and
$(U_yf)(t)=t\,f'(t)$, that is, $U_y=t\frac{d}{dt}$.

\vskip2mm
\textbf{Example 5.}
Let the operator $A(x,y)$ be defined by
$$
(A(x,y)f)(t)=e^{xt+y^2}f(t).
$$
Set $U(x,y)=A(x,y)$.
Then $(U_xf)(t)=t\,e^{xt+y^2}f(t)$, that is,
$U_x=M_{\,t e^{xt+y^2}}$,
where $M_g$ denotes the operator of multiplication by the function $g$.
Moreover,
$(U_yf)(t)=2y\,e^{xt+y^2}f(t)$, that is,
$U_y=M_{\,2y e^{xt+y^2}}$.

\vskip2mm
\textbf{Example 6.}
Let
$$
(U(x,y)f)(t)= \int_0^t e^{x(t-\tau)+y\tau}f(\tau)\,d\tau.
$$
Then
$$
(U_xf)(t)=\int_0^t (t-\tau)e^{x(t-\tau)+y\tau} f(\tau)\,d\tau,
$$
and
$$
(U_yf)(t)=\int_0^t \tau\,e^{x(t-\tau)+y\tau} f(\tau)\,d\tau.
$$

\vskip2mm
\textbf{Example 7.}
Let the Euler operator with variable parameter be defined by
$$
(E_\alpha f)(t)=t^\alpha f'(t),
$$
where
$\alpha=\alpha(x,y)$.
Set
$U(x,y)=E_{\alpha(x,y)}$.

Then
$$
(U_xf)(t)=\alpha_x(x,y)\,\ln t\,t^{\alpha(x,y)}f'(t),
$$
that is,
$$
U_x= \alpha_x(x,y)\,(\ln t)\,t^{\alpha(x,y)}\frac{d}{dt}.
$$
Similarly,
$$
(U_yf)(t)=\alpha_y(x,y)\,\ln t\,t^{\alpha(x,y)}f'(t),
$$
that is,
$$
U_y=\alpha_y(x,y)\,(\ln t)\,t^{\alpha(x,y)}\frac{d}{dt}.
$$

\vskip2mm
\textbf{Example 8.}
Let $X=C^\infty(\mathbb R^2)$,
and let the operator-valued function $U(x,y)$ be defined by
$$
(U(x,y)f)(t_1,t_2)=e^{xt_1+yt_2}\left(\frac{\partial f}{\partial t_1}(t_1,t_2)
+\frac{\partial f}{\partial t_2}(t_1,t_2)\right).
$$
Then
$$
(U_xf)(t_1,t_2)=t_1 e^{xt_1+yt_2} \left( \frac{\partial f}{\partial t_1}
+\frac{\partial f}{\partial t_2}\right)(t_1,t_2),
$$
that is,
$$
U_x=t_1e^{xt_1+yt_2}\left(\frac{\partial}{\partial t_1}
+\frac{\partial}{\partial t_2}\right),
$$
and
$$
(U_yf)(t_1,t_2)=t_2 e^{xt_1+yt_2}\left(\frac{\partial f}{\partial t_1}
+\frac{\partial f}{\partial t_2}\right)(t_1,t_2),
$$
that is,
$$
U_y=t_2e^{xt_1+yt_2}\left(\frac{\partial}{\partial t_1}
+\frac{\partial}{\partial t_2}\right).
$$

\vskip2mm
\textbf{Example 9.}
Let $X=C^\infty(\Omega)$,
where $\Omega=\{(t,\tau)\in\mathbb R^2\colon t\ge0,\ \tau\ge0\}$,
and let the operator-valued function $U(x,y)$ be defined by
$$
(U(x,y)f)(t,\tau)=\int_0^t\int_0^\tau e^{x(t-s)+y(\tau-\sigma)}f(s,\sigma)\,d\sigma\,ds.
$$
Then
$$
(U_xf)(t,\tau)=\int_0^t\int_0^\tau(t-s)\,e^{x(t-s)+y(\tau-\sigma)}f(s,\sigma)\,d\sigma\,ds,
$$
and
$$
(U_yf)(t,\tau)=\int_0^t\int_0^\tau(\tau-\sigma)\,e^{x(t-s)+y(\tau-\sigma)}f(s,\sigma)\,d\sigma\,ds.
$$

Thus, the partial derivatives $U_x$ and $U_y$ are also integral operators.

\vskip2mm
\textbf{Example 10.}
Let
$$
(T_{x,y}f)(t,\tau)=f(t+x,\tau+y).
$$
Then $U(x,y)=T_{x,y}$.
We have
$$
(U_xf)(t,\tau)=\frac{\partial f}{\partial t}(t+x,\tau+y),
$$
that is,
$$
U_x=T_{x,y}\frac{\partial}{\partial t},
$$
and similarly,
$$
U_y=T_{x,y}\frac{\partial}{\partial \tau}.
$$

\vskip2mm
\textbf{Remark 3.}
In all the examples above, the partial derivatives of operator-valued functions are computed by differentiating the dependence of the operators on the external variables $x,y$.
The internal variables $t_1,t_2,\ldots,t_r$ on which the operators act are treated as parameters and remain unchanged under differentiation with respect to $x,y$.

\vskip2mm
\textbf{Remark 4.}
In what follows, continuity of operator-valued functions will be understood in terms of their action on arbitrary elements of the space $X$.
That is, an operator-valued function $U:D\to\mathcal{A}(X)$ belongs to the class
$C^k(D,\mathcal{A}(X))$ if, for every $f\in X$, the function
$u_f(x,y):=(U(x,y)f)$ belongs to the class $C^k(D,X)$.

\vskip2mm
\textbf{Definition 4.}
Let the components $U,V$ of the function (\ref{rozklad za basysom -2-vym-algebra-zi-zminnymy-operatoramy})
have continuous first-order partial derivatives in the domain $D$. The function $\Phi$ will be called monogenic in $D$
if, at every point $(x,y)\in D$, there exists an element
$\Phi'(Z)\in\mathfrak{A}_2(\mathcal{A}(X),\star;x,y)$ such that
\begin{equation}\label{ozn-monog-A-2-zmin-operatory}
d\Phi=\Phi'(Z)dZ.
\end{equation}
The element $\Phi'(Z)$ will be called the derivative of the function $\Phi$ at the point $Z$.

We have
$$
d\Phi=\Phi_xdx+\Phi_ydy, \qquad
dZ=dx E_1+dy E_2.
$$
Therefore, from equality (\ref{ozn-monog-A-2-zmin-operatory}), we obtain
$$
\Phi_xdx+\Phi_ydy =\Phi'(Z)dx E_1+\Phi'(Z)dy E_2.
$$
Since $E_1$ is the identity element of the algebra, comparison of the coefficients of $dx$ gives
$\Phi'(Z)=\Phi_x$.
Substituting this expression into the equality of the coefficients of $dy$, we obtain
$\Phi_y=\Phi_xE_2$.

Thus, we have the following theorem.

\vskip2mm
\textbf{Theorem 1.}
A function $\Phi$ of the form (\ref{rozklad za basysom -2-vym-algebra-zi-zminnymy-operatoramy}) is monogenic in the domain $D$ if and only if
\begin{equation}\label{A-2-zmin-operatory-umovy-K-R}
\Phi_y=\Phi_xE_2.
\end{equation}
Conditions (\ref{A-2-zmin-operatory-umovy-K-R}) are analogues of the Cauchy--Riemann conditions.

Let us write the Cauchy--Riemann conditions in terms of the components $U,V$ of a monogenic function.

We have
$$
\Phi_x=U_xE_1+V_xE_2,
$$
and therefore, from condition (\ref{A-2-zmin-operatory-umovy-K-R}) and the multiplication table of the algebra, we obtain
$$
U_yE_1+V_yE_2
=
A(x,y)V_xE_1+\bigl(U_x+B(x,y)V_x\bigr)E_2.
$$
Hence, equating the coefficients of $E_1$ and $E_2$, we obtain the system
\begin{equation}\label{A-2-zmin-operatopy-systema-K-R}
U_y=A(x,y)V_x,\qquad V_y=U_x+B(x,y)V_x.
\end{equation}

Clearly, the Cauchy--Riemann conditions (\ref{A-2-zmin-operatory-umovy-K-R}) or (\ref{A-2-zmin-operatopy-systema-K-R}) depend on the point $(x,y)\in D$.
Thus, in general, the Cauchy--Riemann conditions are different at different points of $D$.

\vskip2mm
\textbf{Example 11.}
Consider the function $\Phi(Z)=E_2^2$.

According to the multiplication table (\ref{tabl-mnozh-alg-2-zmin-operat}) of the algebra
$\mathfrak{A}_2(\mathcal{A}(X),\star;x,y)$, we have
$$
E_2^2=A(x,y)E_1+B(x,y)E_2.
$$
Hence,
$\Phi(Z)=A(x,y)E_1+B(x,y)E_2$, that is,
$$
U(x,y)=A(x,y),\qquad V(x,y)=B(x,y).
$$

Substituting these functions into the Cauchy--Riemann conditions
(\ref{A-2-zmin-operatopy-systema-K-R}), we obtain the system
\begin{equation}\label{A-2-zmin-operatory-umova-symisnosti}
A_y=A\,B_x, \qquad B_y=A_x+B\,B_x.
\end{equation}

Thus, the function $\Phi(Z)=E_2^2$ is monogenic in the domain $D$ if and only if the structural operators
$A(x,y)$ and $B(x,y)$ possess all first-order partial derivatives and satisfy system (\ref{A-2-zmin-operatory-umova-symisnosti}).

Thus, even such a simple function as $E_2^2$ is not always monogenic.
For its monogenicity, it is necessary and sufficient that the structural operators of the algebra satisfy system
(\ref{A-2-zmin-operatory-umova-symisnosti}).

\vskip2mm
\textbf{Definition 5.}
System (\ref{A-2-zmin-operatory-umova-symisnosti}) will be called the \textbf{compatibility system}.
The algebra $\mathfrak{A}_2(\mathcal{A}(X),\star;x,y)$ will be called \textbf{compatible} if its structural operators
$A(x,y)$ and $B(x,y)$ possess all first-order partial derivatives and satisfy the compatibility system
(\ref{A-2-zmin-operatory-umova-symisnosti}).

In what follows, we consider only compatible algebras.

Let us investigate some properties of monogenic functions.

\vskip2mm
\vskip2mm
\textbf{Theorem 2.}
Let the functions $\Phi$ and $\Psi$ be monogenic in the domain $D$ in a compatible algebra
$\mathfrak{A}_2(\mathcal{A}(X),\star;x,y)$.
Then, for arbitrary constant operators $\alpha,\beta\in\mathcal{A}(X)$, the function
$$
\Theta=\alpha\Phi+\beta\Psi
$$
is also monogenic in the domain $D$.

\textbf{Proof.}
Since the functions $\Phi$ and $\Psi$ are monogenic in $D$, they satisfy the Cauchy--Riemann conditions
$$
\Phi_y=\Phi_xE_2,\qquad \Psi_y=\Psi_xE_2.
$$
Consider the function
$$
\Theta=\alpha\Phi+\beta\Psi,\qquad
\alpha,\beta\in\mathcal{A}(X),
$$
where the operators $\alpha$ and $\beta$ do not depend on $(x,y)$. Then
$$
\Theta_y=(\alpha\Phi+\beta\Psi)_y
=\alpha\Phi_y+\beta\Psi_y.
$$

Using the monogenicity of the functions $\Phi$ and $\Psi$, we obtain
$$
\Theta_y=\alpha\Phi_xE_2+\beta\Psi_xE_2.
$$
Since $\mathcal{A}(X)$ is a commutative associative algebra of operators,
$$
\alpha\Phi_xE_2+\beta\Psi_xE_2
=(\alpha\Phi_x+\beta\Psi_x)E_2.
$$
On the other hand,
$$
\Theta_x=(\alpha\Phi+\beta\Psi)_x
=\alpha\Phi_x+\beta\Psi_x.
$$
Hence,
$\Theta_y=\Theta_xE_2$.

Therefore, the function $\Theta$ is monogenic in the domain $D$.
The theorem is proved.

\vskip2mm
\textbf{Theorem 3.}
Let the functions $\Phi$ and $\Psi$ be monogenic in the domain $D$ in a compatible algebra
$\mathfrak{A}_2(\mathcal{A}(X),\star;x,y)$.
Then their product $\Theta=\Phi\star\Psi$ is also monogenic in $D$. Moreover,
$$
(\Phi\star\Psi)'=\Phi'\star\Psi+\Phi\star\Psi'.
$$

\textbf{Proof.}
Let
$$
\Phi=U E_1+V E_2,\qquad
\Psi=P E_1+Q E_2,
$$
where $U,V,P,Q\in C^1(D,\mathcal{A}(X))$.

Since $E_2^2=A(x,y)E_1+B(x,y)E_2$, we have
$$
\Theta=\Phi\star\Psi=W E_1+R E_2,
$$
where
$$
W=UP+AVQ,\qquad R=UQ+VP+BVQ.
$$

Since the functions $\Phi$ and $\Psi$ are monogenic, their components satisfy the Cauchy--Riemann conditions:
$$
U_y=A V_x,\qquad V_y=U_x+B V_x,
$$
$$
P_y=A Q_x,\qquad Q_y=P_x+B Q_x.
$$

We show that the function $\Theta$ also satisfies the Cauchy--Riemann conditions, namely,
$$
W_y=A R_x,\qquad R_y=W_x+B R_x.
$$
We have
$$
W_y=U_yP+UP_y+A_yVQ+A V_yQ+A VQ_y.
$$
Substituting the Cauchy--Riemann conditions for $\Phi$ and $\Psi$, we obtain
$$
W_y=A V_xP+UAQ_x+A_yVQ+A(U_x+B V_x)Q+AV(P_x+B Q_x).
$$
On the other hand,
$$
R_x=U_xQ+UQ_x+V_xP+VP_x+B_xVQ+B V_xQ+B VQ_x.
$$
Therefore,
$$
A R_x=AU_xQ+AUQ_x+AV_xP+AVP_x+AB_xVQ+ABV_xQ+ABVQ_x.
$$

Since the algebra $\mathcal{A}(X)$ is commutative, comparison of the two expressions yields
$$
W_y-AR_x=(A_y-AB_x)VQ.
$$
In a compatible algebra, the structural operators satisfy the compatibility condition
$A_y=AB_x$.

Hence, $W_y=AR_x$.

Now let us verify the second Cauchy--Riemann condition. We have
$$
R_y=U_yQ+UQ_y+V_yP+VP_y+B_yVQ+B V_yQ+B VQ_y.
$$
Substituting the Cauchy--Riemann conditions, we obtain
$$
R_y=AV_xQ+U(P_x+BQ_x)+(U_x+BV_x)P+VAQ_x
$$
$$
+B_yVQ+B(U_x+BV_x)Q+BV(P_x+BQ_x).
$$
On the other hand,
$$
W_x=U_xP+UP_x+A_xVQ+AV_xQ+AVQ_x,
$$
and
$$
BR_x=BU_xQ+BUQ_x+BV_xP+BVP_x+BB_xVQ+B^2V_xQ+B^2VQ_x.
$$

Therefore, $W_x+BR_x$ is equal to
$$
U_xP+UP_x+A_xVQ+AV_xQ+AVQ_x+BU_xQ+BUQ_x+BV_xP
$$
$$
+BVP_x+BB_xVQ+B^2V_xQ+B^2VQ_x.
$$
Comparing this expression with the formula for $R_y$, we obtain
$$
R_y-(W_x+BR_x)=(B_y-A_x-BB_x)VQ.
$$
Since, in a compatible algebra,
$B_y=A_x+BB_x$, it follows that
$$
R_y=W_x+BR_x.
$$

Thus, the components $W$ and $R$ of the function
$\Theta=\Phi\star\Psi$ satisfy the Cauchy--Riemann system. Therefore, the function $\Theta$ is monogenic in $D$.

It remains to prove the derivative formula. Since
$$
\Phi'=\Phi_x,\qquad \Psi'=\Psi_x,
$$
the product rule gives
$$
(\Phi\star\Psi)'=(\Phi\star\Psi)_x
=\Phi_x\star\Psi+\Phi\star\Psi_x
=\Phi'\star\Psi+\Phi\star\Psi'.
$$
The theorem is proved.

\vskip2mm
\textbf{Remark 5.}
Theorem 3 shows that the compatibility system (\ref{A-2-zmin-operatory-umova-symisnosti}) is a sufficient condition for the set of monogenic functions to be closed under $\star$-multiplication.
In general, however, condition (\ref{A-2-zmin-operatory-umova-symisnosti}) is not necessary.
We note that, in the general operator case, the converse statement requires additional assumptions on the algebra $\mathcal{A}(X)$.
One possible assumption is that the commutative associative algebra of operators $\mathcal{A}(X)$ contains no zero divisors.
The investigation of such conditions, however, lies beyond the scope of the present paper.

Let us consider some further examples of monogenic functions.

\vskip2mm
\textbf{Example 12.}
The function $\Phi(Z)=Z=xE_1+yE_2$ is monogenic in the domain $D$.

Indeed,
$\Phi_x=E_1$, $\Phi_y=E_2$.
Since $E_1$ is the identity element of the algebra,
$$
\Phi_xE_2=E_1E_2=E_2=\Phi_y.
$$

Hence, the Cauchy--Riemann condition
$\Phi_y=\Phi_xE_2$ is satisfied.
Therefore, the function $\Phi(Z)=Z$ is monogenic in $D$. Moreover,
$\Phi'(Z)=\Phi_x=E_1$.

\vskip2mm
\textbf{Example 13.}\label{pryklad-Z-star-2}
Consider the function $\Phi(Z)=Z^{\star2}$.

Since the function $Z$ is monogenic and the product of two monogenic functions in a compatible algebra is again monogenic,
$Z^{\star2}=Z\star Z$ is monogenic in $D$.

Moreover,
$Z^{\star2}=(xE_1+yE_2)\star(xE_1+yE_2)$, and therefore
$$
Z^{\star2}
=
\bigl(x^2I+y^2A(x,y)\bigr)E_1
+
\bigl(2xyI+y^2B(x,y)\bigr)E_2.
$$

The following corollary follows from the preceding two theorems.

\vskip2mm
\textbf{Corollary 1.}\label{naslidok-1-A-2-oper}
The set of all functions monogenic in the domain $D$ with values in the compatible algebra
$\mathfrak{A}_2(\mathcal A(X),\star;x,y)$
forms a commutative associative algebra with respect to addition of functions and $\star$-multiplication.

Corollary \ref{naslidok-1-A-2-oper} provides a method for constructing monogenic functions.

It follows from Example \ref{pryklad-Z-star-2} and Theorem 3 that, for every natural number $n$, the function
$Z^{\star n}$ is monogenic in the domain $D$.

Theorems 2 and 3 now imply that, for arbitrary constant operators
$C_0,C_1,\ldots,C_n\in\mathcal{A}(X)$,
the polynomial
$$
P_n(Z)=\sum_{k=0}^{n} C_kZ^{\star k}
$$
is monogenic in the domain $D$.
In particular, if $C_k=c_kI$, $c_k\in\mathbb{R}$, we obtain the scalar polynomial
$$
P_n(Z)=\sum_{k=0}^{n} c_k Z^{\star k},
$$
which is also monogenic.

\vskip2mm
\textbf{Theorem 4.}
Let the structural constants $A$ and $B$ of the compatible algebra
$\mathfrak{A}_2(\mathcal{A}(X),\star;x,y)$ be constant.
If the function
$\Phi(x,y)=U(x,y)E_1+V(x,y)E_2$
is monogenic and its components $U,V$ possess continuous partial derivatives up to and including second order, then its derivative
$$
\Phi'(Z)=\Phi_x(x,y)=U_x(x,y)E_1+V_x(x,y)E_2
$$
is also monogenic.

\textbf{Proof.}
Since $\Phi$ is monogenic, its components satisfy the Cauchy--Riemann conditions
(\ref{A-2-zmin-operatopy-systema-K-R}):
$$
U_y=AV_x,\qquad V_y=U_x+BV_x.
$$
Differentiating these equalities with respect to $x$, we obtain
$$
U_{xy}=AV_{xx}, \qquad
V_{xy}=U_{xx}+BV_{xx}.
$$
These are precisely the monogenicity conditions for the function
$\Phi_x=U_xE_1+V_xE_2$.
The theorem is proved.

\vskip2mm
\textbf{Remark 6.}
If the structural constants $A$ and $B$ are operator-valued functions of $x,y$, then the preceding theorem is, in general, no longer valid.

\section{Relationship Between Monogenic Functions and Operator-Differential Equations}\label{Paragraf-3-5}

In this section, we establish a relationship between monogenic functions in the compatible algebra
$\mathfrak{A}_2(\mathcal{A}(X),\star;x,y)$ and operator-differential equations with variable coefficients. This makes it possible to construct explicit solutions of broad classes of operator-differential, integro-differential, and other functional equations by means of the theory of monogenic functions.

\vskip2mm
\textbf{Theorem 5.}\label{zvyazok-z-rivnianniam-zmin-operatory}
Let $\Phi(Z)=U(x,y)E_1+V(x,y)E_2$ be a monogenic function in a domain $D$ with values in the compatible algebra $\mathfrak{A}_2(\mathcal{A}(X),\star;x,y)$,
and let $U,V\in C^2(D,\mathcal{A}(X))$. Then the component $V$ satisfies the operator-differential equation
\begin{equation}\label{rivn-V-zmin-operatory}
V_{yy}=A V_{xx}+B V_{xy}+(2A_x+B B_x)V_x,
\end{equation}
whereas the component $U$ satisfies the operator-differential equation
\begin{equation}\label{rivn-U-zmin-operatory}
A U_{yy}=A^2 U_{xx}+A B U_{xy}+(2A B_x-B A_x)U_y.
\end{equation}

\textbf{Proof.} Since the function $\Phi$ is monogenic, its components satisfy the Cauchy--Riemann system:
\begin{equation}\label{A-2-zmin-operatopy-systema-K-R-vyvid}
U_y=AV_x,\qquad V_y=U_x+BV_x.
\end{equation}
First, we derive the equation for the component $V$. Differentiating the second equality in
(\ref{A-2-zmin-operatopy-systema-K-R-vyvid}) with respect to $y$, we obtain
$$
V_{yy}=U_{xy}+B_yV_x+BV_{xy}.
$$
Differentiating the first equality in (\ref{A-2-zmin-operatopy-systema-K-R-vyvid}) with respect to $x$, we have
$$
U_{xy}=(AV_x)_x=A_xV_x+AV_{xx}.
$$
Hence,
$$
V_{yy}=AV_{xx}+BV_{xy}+(A_x+B_y)V_x.
$$
Since the algebra is compatible, we have
$B_y=A_x+B B_x$.

Therefore,
$$
V_{yy}=A V_{xx}+B V_{xy}+(2A_x+B B_x)V_x.
$$

Now we derive the equation for the component $U$. Differentiating the first equality in
(\ref{A-2-zmin-operatopy-systema-K-R-vyvid}) with respect to $y$, we obtain
$$
U_{yy}=A_yV_x+AV_{xy}.
$$
Differentiating the second equality in (\ref{A-2-zmin-operatopy-systema-K-R-vyvid}) with respect to $x$, we have
$$
V_{xy}=U_{xx}+B_xV_x+BV_{xx}.
$$
Therefore,
$$
U_{yy}=A_yV_x+AU_{xx}+AB_xV_x+ABV_{xx}.
$$
On the other hand, from the equality
$U_{xy}=A_xV_x+AV_{xx}$ we obtain
$$
AV_{xx}=U_{xy}-A_xV_x.
$$
Hence,
$$
ABV_{xx}=BU_{xy}-BA_xV_x.
$$
Therefore,
$$
U_{yy}=AU_{xx}+BU_{xy}+(A_y+AB_x-BA_x)V_x.
$$
Since the algebra is compatible,
$A_y=AB_x$.
Thus,
$$
U_{yy}=AU_{xx}+BU_{xy}+(2AB_x-BA_x)V_x.
$$
Applying the operator $A$ to both sides of this equality, we obtain
$$
A U_{yy}=A^2U_{xx}+ABU_{xy}+(2AB_x-BA_x)AV_x.
$$
But from the first equation of the Cauchy--Riemann system, we have
$AV_x=U_y$.
Therefore,
$$
A U_{yy}=A^2 U_{xx}+A B U_{xy}+(2A B_x-B A_x)U_y.
$$
The theorem is proved.

We now present converse results.

\vskip2mm
\textbf{Theorem 6.}\label{obernena-teorema-U-zmin-operatory}
Let $D\subset\mathbb{R}^2$ be a simply connected domain, and let the operator-valued functions
$A,B\in C^1(D,\mathcal{A}(X))$
define a compatible algebra $\mathfrak{A}_2(\mathcal{A}(X),\star;x,y)$. Suppose, in addition, that for every point $(x,y)\in D$, the operator $A(x,y)$
is invertible in the algebra $\mathcal{A}(X)$, and let the operator-valued function $U\in C^2(D,\mathcal{A}(X))$ satisfy the equation
\begin{equation}\label{U-zmin-operatory}
A U_{yy}=A^2 U_{xx}+AB U_{xy}+(2AB_x-BA_x)U_y.
\end{equation}
Then there exists an operator-valued function $V=V(x,y)$ such that
$\Phi(Z)=U(x,y)E_1+V(x,y)E_2$
is monogenic in the domain $D$.

\textbf{Proof.}
For the function $\Phi=UE_1+VE_2$ to be monogenic, its components must satisfy the Cauchy--Riemann system
(\ref{A-2-zmin-operatopy-systema-K-R-vyvid}).
Since the operator $A$ is invertible, the first Cauchy--Riemann equation gives
$$
V_x=A^{-1}U_y.
$$
Then the second equation takes the form
$$
V_y=U_x+BA^{-1}U_y.
$$
Thus, the function $V$ must be found from the system
$$
V_x=A^{-1}U_y,\qquad
V_y=U_x+BA^{-1}U_y.
$$

Consider the operator-valued differential form
$$
\omega=A^{-1}U_y\,dx+\left(U_x+BA^{-1}U_y\right)dy.
$$
The function $V$ exists if the form $\omega$ is exact. Since the domain $D$ is simply connected, it is sufficient to verify that $\omega$ is closed:
$$
(A^{-1}U_y)_y=\left(U_x+BA^{-1}U_y\right)_x.
$$

Let us compute the left-hand side. Since
$(A^{-1})_y=-A^{-1}A_yA^{-1}$, we have
$$
(A^{-1}U_y)_y=A^{-1}U_{yy}-A^{-1}A_yA^{-1}U_y.
$$
The right-hand side is
$$
\left(U_x+BA^{-1}U_y\right)_x
=
U_{xx}+B_xA^{-1}U_y+B(A^{-1})_xU_y+BA^{-1}U_{xy}.
$$
Since
$(A^{-1})_x=-A^{-1}A_xA^{-1}$, we obtain
$$
\left(U_x+BA^{-1}U_y\right)_x
=
U_{xx}+B_xA^{-1}U_y-BA^{-1}A_xA^{-1}U_y+BA^{-1}U_{xy}.
$$

Hence, the closedness condition takes the form
$$
A^{-1}U_{yy}- A^{-1}A_yA^{-1}U_y
=
U_{xx}+B_xA^{-1}U_y-BA^{-1}A_xA^{-1}U_y+BA^{-1}U_{xy}.
$$
Multiplying this equality on the left by $A^2$, we obtain
$$
AU_{yy}-A_yU_y
=
A^2U_{xx}+AB_xU_y-BA_xU_y+ABU_{xy}.
$$
Hence,
$$
AU_{yy}
=
A^2U_{xx}+ABU_{xy}+(A_y+AB_x-BA_x)U_y.
$$
Since the algebra is compatible, $A_y=AB_x$.
Therefore,
$$
AU_{yy}
=
A^2U_{xx}+ABU_{xy}+(2AB_x-BA_x)U_y.
$$
But this is exactly equation (\ref{U-zmin-operatory}), which is satisfied by $U$ by assumption.

Hence, the form $\omega$ is closed. Therefore, in the simply connected domain $D$, there exists an operator-valued function $V$ such that
$dV=\omega$. Thus,
$$
V_x=A^{-1}U_y,\qquad
V_y=U_x+BA^{-1}U_y.
$$
Therefore,
$$
U_y=AV_x,\qquad
V_y=U_x+BV_x.
$$
Hence, the function
$\Phi(Z)=U(x,y)E_1+V(x,y)E_2$
is monogenic in the domain $D$.
The theorem is proved.

Moreover, the function $V$ can be represented by the formula
$$
V(x,y)=V(x_0,y_0)+\int_\gamma A^{-1}U_y\,dx+\left(U_x+BA^{-1}U_y\right)dy,
$$
where $\gamma$ is an arbitrary curve in $D$ joining the point $(x_0,y_0)$ to the point $(x,y)$. Since the form is closed, this integral is independent of the choice of path.

\vskip2mm
\textbf{Theorem 7.}\label{obernena-teorema-V-zmin-operatory}
Let $D\subset\mathbb{R}^2$ be a simply connected domain, and let the operator-valued functions
$A,B\in C^1(D,\mathcal{A}(X))$
define a compatible algebra $\mathfrak{A}_2(\mathcal{A}(X),\star;x,y)$.
Suppose, in addition, that
$V\in C^2(D,\mathcal{A}(X))$ satisfies the equation
\begin{equation}\label{V-zmin-operatory}
V_{yy}=A V_{xx}+B V_{xy}+(2A_x+B B_x)V_x.
\end{equation}
Then there exists an operator-valued function $U=U(x,y)$ such that
$\Phi(Z)=U(x,y)E_1+V(x,y)E_2$
is monogenic in the domain $D$.

\textbf{Proof.}
For the function $\Phi=UE_1+VE_2$ to be monogenic, its components must satisfy the Cauchy--Riemann system
(\ref{A-2-zmin-operatopy-systema-K-R-vyvid}).
From the second Cauchy--Riemann equation, we obtain
$$
U_x=V_y-BV_x.
$$
Thus, the function $U$ must be found from the system
$$
U_x=V_y-BV_x,\qquad U_y=AV_x.
$$

Consider the operator-valued differential form
$$
\omega=(V_y-BV_x)\,dx+AV_x\,dy.
$$
The function $U$ exists if the form $\omega$ is exact. Since the domain $D$ is simply connected, it is sufficient to verify that $\omega$ is closed:
$$
(V_y-BV_x)_y=(AV_x)_x.
$$

We compute
$$
(V_y-BV_x)_y=V_{yy}-B_yV_x-BV_{xy},
$$
and
$$
(AV_x)_x=A_xV_x+AV_{xx}.
$$
Therefore, the closedness condition takes the form
$$
V_{yy}-B_yV_x-BV_{xy}
=
A_xV_x+AV_{xx}.
$$
Hence,
$$
V_{yy}=AV_{xx}+BV_{xy}+(A_x+B_y)V_x.
$$
Since the algebra is compatible,
$B_y=A_x+BB_x$.
Therefore,
$$
V_{yy}=A V_{xx}+B V_{xy}+(2A_x+B B_x)V_x.
$$
But this is precisely the equation satisfied by the function $V$ by assumption.

Hence, the form $\omega$ is closed. Therefore, in the simply connected domain $D$, there exists an operator-valued function $U$ such that
$dU=\omega$.
That is,
$$
U_x=V_y-BV_x,\qquad U_y=AV_x.
$$
Hence,
$$
U_y=AV_x,\qquad V_y=U_x+BV_x.
$$
Therefore, the function
$\Phi(Z)=U(x,y)E_1+V(x,y)E_2$
is monogenic in the domain $D$.
The theorem is proved.

Moreover, the function $U$ can be represented by the formula
$$
U(x,y)=U(x_0,y_0)+\int_\gamma (V_y-BV_x)\,dx+AV_x\,dy,
$$
where $\gamma$ is an arbitrary curve in $D$ joining the point $(x_0,y_0)$ to the point $(x,y)$. Since the form is closed, this integral is independent of the choice of path.

\vskip2mm
\textbf{Remark 7.}
Thus, the set of twice continuously differentiable operator-valued solutions of equations
(\ref{rivn-V-zmin-operatory}) and (\ref{rivn-U-zmin-operatory})
in a simply connected domain $D$ is in one-to-one correspondence with the set of monogenic functions in the compatible algebra
$\mathfrak{A}_2(\mathcal{A}(X),\star;x,y)$.

More precisely, each monogenic function
$\Phi(Z)=U(x,y)E_1+V(x,y)E_2$
corresponds to a pair of operator-valued functions $U$ and $V$ satisfying equations
(\ref{rivn-U-zmin-operatory}) and (\ref{rivn-V-zmin-operatory}), respectively.

Conversely, for every solution
$U\in C^2(D,\mathcal{A}(X))$
of equation (\ref{rivn-U-zmin-operatory}), there exists a function $V$, determined up to an additive constant, such that
$\Phi(Z)=U(x,y)E_1+V(x,y)E_2$
is monogenic in the domain $D$.

Similarly, for every solution
$V\in C^2(D,\mathcal{A}(X))$
of equation (\ref{rivn-V-zmin-operatory}), there exists a function $U$, determined up to an additive constant, such that
$\Phi(Z)=U(x,y)E_1+V(x,y)E_2$
is monogenic in the domain $D$.

Thus, all solutions of equations
(\ref{rivn-V-zmin-operatory}) and (\ref{rivn-U-zmin-operatory})
in the class $C^2(D,\mathcal{A}(X))$
can be obtained as components of monogenic functions with values in the compatible algebra
$\mathfrak{A}_2(\mathcal{A}(X),\star;x,y)$.

Theorems 6 and 7
provide a constructive method for recovering a monogenic function from a prescribed component.
\vskip2mm

\textbf{Remark.} Second-order operator-differential equations have been studied in numerous works. For methods of solving such equations, see, for example, \cite{Keyantuo,Pyatkov,Gorbachuk}.
\vskip2mm

\section{Examples}\label{Paragraf-6-zmin-operatory}

We give several examples of operator-valued functions $A$ and $B$ satisfying the compatibility system (\ref{A-2-zmin-operatory-umova-symisnosti})
and the corresponding operator-differential equations for the components of monogenic functions.

Equations (\ref{rivn-V-zmin-operatory}) and (\ref{rivn-U-zmin-operatory}) are linear operator-differential equations in the class of
operator-valued functions. However, by realizing the operator-valued functions $A(x,y)$ and $B(x,y)$ on specific linear topological spaces $X$,
we obtain different classes of differential, integro-differential, and functional-differential equations or systems of equations of different
orders and involving different numbers of variables for $\mathbb{K}$-valued functions. We illustrate this by concrete examples.

As noted above, the variables $x$ and $y$ denote the coordinates of a point in the domain $D\subset\mathbb{R}^2$, whereas the variables
$t_1,t_2$, $s$, $t$, $\tau$, $\eta$, and others are used to describe the space $X$ on which the operators of the algebra $\mathcal{A}(X)$ act.
The operators $A(x,y)$ and $B(x,y)$ may depend both on the variables $x,y$ and on the internal variables of the space $X$.
When differentiating with respect to $x$ and $y$, the internal variables of the space $X$ are regarded as independent of $x,y$.

\vskip2mm
\textbf{Principle of transition from operator-differential equations to $\mathbb{K}$-valued equations.}
The passage from the operator-differential equations (\ref{rivn-V-zmin-operatory}), (\ref{rivn-U-zmin-operatory}) to $\mathbb{K}$-valued equations is carried out by applying both sides
of the operator equation to an arbitrary element $f\in X$. One then evaluates both sides of the resulting equality with respect to the internal variables of the space $X$.

For example, if the operator-valued function $V$ satisfies the equation
$$
V_{yy}=AV_{xx}+BV_{xy}+(2A_x+BB_x)V_x,
$$
then, applying both sides of this equation to an arbitrary element
$$
f=f(s,\tau,\ldots)\in X,
$$
we obtain
$$
V_{yy}f=AV_{xx}f+BV_{xy}f+(2A_x+BB_x)V_xf.
$$
We then define the $\mathbb{K}$-valued function
$$
v(x,y,s,\tau,\ldots)=(V(x,y)f)(s,\tau,\ldots),
$$
and, by evaluating both sides of the equality with respect to the internal variables of the space $X$, obtain the corresponding $\mathbb{K}$-valued differential equation.

\vskip2mm
\textbf{Example 13.} Let $X=C^\infty(\mathbb{R})$, and let $I$ be the identity operator. Set
$$
A(x,y)=-\frac{x}{1-y}I,\qquad
B(x,y)=\left(1+\frac{x}{1-y}\right)I,\qquad y\neq1.
$$
Then
$$
A_x=-\frac{1}{1-y}I,\qquad
B_x=\frac{1}{1-y}I.
$$
A direct verification shows that $A$ and $B$ satisfy the compatibility system (\ref{A-2-zmin-operatory-umova-symisnosti}).
Hence, the corresponding algebra is compatible.

The \textbf{operator equation} for the component $U(x,y)$ has the form
$$
-\frac{x}{1-y}\,U_{yy}I
=
\frac{x^2}{(1-y)^2}\,U_{xx}I
-\frac{x(1-y+x)}{(1-y)^2}\,U_{xy}I
+\frac{1-y-x}{(1-y)^2}\,U_yI.
$$

The \textbf{operator equation} for the component $V(x,y)$ has the form
$$
V_{yy}I
=
-\frac{x}{1-y}\,V_{xx}I
+\left(1+\frac{x}{1-y}\right)V_{xy}I
+\frac{x+y-1}{(1-y)^2}\,V_xI.
$$

Since $I$ is the identity operator, according to the principle of transition from
operator-differential equations to $\mathbb K$-valued equations, we apply
both sides of each of the above equations to an arbitrary element
$f\in X$ and define
$$
u(x,y,s)=(U(x,y)f)(s),\qquad
v(x,y,s)=(V(x,y)f)(s).
$$
Since
$$
(Ig)(s)=g(s)
$$
for every function $g\in X$, the action of the identity operator does not change
the form of the equations. Thus, we obtain the following $\mathbb K$-valued
differential equations.

The \textbf{$\mathbb{K}$-valued equation} for the component $u(x,y,s)$ has the form
$$
-\frac{x}{1-y}\,u_{yy}
=
\frac{x^2}{(1-y)^2}\,u_{xx}
-\frac{x(1-y+x)}{(1-y)^2}\,u_{xy}
+\frac{1-y-x}{(1-y)^2}\,u_y.
$$

The \textbf{$\mathbb{K}$-valued equation} for the component $v(x,y,s)$ has the form
$$
v_{yy}
=
-\frac{x}{1-y}\,v_{xx}
+\left(1+\frac{x}{1-y}\right)v_{xy}
+\frac{x+y-1}{(1-y)^2}\,v_x.
$$

\vskip2mm
\textbf{Example 14.} Let $X=C^\infty(\mathbb{R})$ be the space of all infinitely differentiable functions $f=f(s)$, and let $M_g$ denote the operator of multiplication by a function $g=g(s)$:
$$
(M_g f)(s)=g(s)f(s).
$$
\vskip2mm

\textbf{1.}
Let $\alpha=\alpha(s)$ be a fixed smooth function. Set
$$
A(x,y)=M_{-\frac{\alpha(s)x}{1-\alpha(s)y}},
\qquad
B(x,y)=M_{1+\frac{\alpha(s)x}{1-\alpha(s)y}}.
$$
Since the product of multiplication operators corresponds to the product of their defining functions, the compatibility system reduces to scalar equalities for each fixed $s$. Hence,
the compatibility system (\ref{A-2-zmin-operatory-umova-symisnosti}) is satisfied.

In this case, the \textbf{operator-differential equation} for the component
$U(x,y)$ has the form
$$
M_{-\frac{\alpha(s)x}{1-\alpha(s)y}}\,U_{yy}
=
M_{\frac{\alpha^2(s)x^2}{(1-\alpha(s)y)^2}}\,U_{xx}
-
M_{\frac{\alpha(s)x\bigl(1-\alpha(s)y+\alpha(s)x\bigr)}
{(1-\alpha(s)y)^2}}\,U_{xy}
+
M_{\frac{\alpha(s)\bigl(1-\alpha(s)y-\alpha(s)x\bigr)}
{(1-\alpha(s)y)^2}}\,U_y.
$$

The \textbf{operator-differential equation} for the component $V(x,y)$ has the form
$$
V_{yy}
=
M_{-\frac{\alpha(s)x}{1-\alpha(s)y}}\,V_{xx}
+
M_{1+\frac{\alpha(s)x}{1-\alpha(s)y}}\,V_{xy}
+
M_{\frac{\alpha(s)\bigl(\alpha(s)x+\alpha(s)y-1\bigr)}
{(1-\alpha(s)y)^2}}\,V_x.
$$

Since $M_{\alpha}$ is the multiplication operator by the function
$\alpha(s)$, according to the principle of transition from
operator-differential equations to $\mathbb K$-valued equations,
we apply both sides of each of the above equations to an arbitrary element
$f\in X$ and define
$$
u(x,y,s)=(U(x,y)f)(s),\qquad
v(x,y,s)=(V(x,y)f)(s).
$$
Here
$$
(M_{\alpha}g)(s)=\alpha(s)g(s)
$$
for every function $g\in X$. Therefore, multiplication operators
are replaced by multiplication by the corresponding functions, and we obtain the following
$\mathbb K$-valued differential equations.

For the component $u(x,y,s)$, we obtain the following \textbf{$\mathbb K$-valued differential equation}:
$$
-\frac{\alpha x}{1-\alpha y}\,u_{yy}
=
\frac{\alpha^2x^2}{(1-\alpha y)^2}\,u_{xx}
-
\frac{\alpha x(1-\alpha y+\alpha x)}{(1-\alpha y)^2}\,u_{xy}
+
\frac{\alpha(1-\alpha y-\alpha x)}{(1-\alpha y)^2}\,u_y.
$$

The \textbf{$\mathbb K$-valued differential equation} for the component $v(x,y,s)$ is
$$
v_{yy}
=
-\frac{\alpha x}{1-\alpha y}\,v_{xx}
+
\left(1+\frac{\alpha x}{1-\alpha y}\right)v_{xy}
+
\frac{\alpha(\alpha x+\alpha y-1)}{(1-\alpha y)^2}\,v_x.
$$
Here $\alpha=\alpha(s)$, and the operators act with respect to the internal variable $s$.

\textbf{2.}
Let $X=C^\infty(\mathbb{R}^2)$ and $\alpha(s,\tau)=e^{-s}\cos\tau$. Then for the component $v(x,y,s,\tau)$ we obtain the equation
$$
v_{yy}
=
-\frac{x e^{-s}\cos\tau}{1-y e^{-s}\cos\tau}\,v_{xx}
+
\left(1+\frac{x e^{-s}\cos\tau}{1-y e^{-s}\cos\tau}\right)\,v_{xy}
$$
$$
+
\frac{e^{-s}\cos\tau\left(
x e^{-s}\cos\tau+y e^{-s}\cos\tau-1
\right)}
{\left(1-y e^{-s}\cos\tau\right)^2}\,v_x.
$$

Here the coefficients depend not only on the independent variables $x,y$, but also on the internal variables $s,\tau$ of the space $X$.

\vskip2mm
\textbf{Example 15.} Let $P$ be a projection operator in the space $X$, that is, $P^2=P$.
For example, in the space $X=C^\infty(\mathbb{R})$, one may take
$$
(Pf)(s)=f(0).
$$
Set
$$
A(x,y)=-\frac{x}{1-y}P,\qquad
B(x,y)=\left(1+\frac{x}{1-y}\right)P.
$$
Since $P^2=P$, the compatibility system again reduces to the scalar system for the coefficients:
$$
A_y=AB_x,\qquad
B_y=A_x+BB_x.
$$

Hence, the \textbf{operator-differential equation} for the component $U(x,y)$ has the form
$$
-\frac{x}{1-y}PU_{yy}
=
\frac{x^2}{(1-y)^2}PU_{xx}
-\frac{x(1-y+x)}{(1-y)^2}PU_{xy}
+\frac{1-y-x}{(1-y)^2}PU_y.
$$

The \textbf{operator-differential equation} for the component $V(x,y)$ has the form
\begin{equation}\label{rivn-pryklad-15-V-dali}
V_{yy}
=
-\frac{x}{1-y}PV_{xx}
+
\left(1+\frac{x}{1-y}\right)PV_{xy}
+
\frac{x+y-1}{(1-y)^2}PV_x.
\end{equation}

According to the principle of transition from operator-differential equations to
$\mathbb K$-valued equations, we apply both sides of each of the above
equations to an arbitrary element $f\in X$ and define
$$
u(x,y,s)=(U(x,y)f)(s),\qquad
v(x,y,s)=(V(x,y)f)(s).
$$
Since $(Pf)(s)=f(0)$, we have
$$
(Pg)(s)=g(0)
$$
for every function $g\in X$. Therefore, the operator $P$ is replaced by
evaluation of the function at the point $s=0$, and we obtain the following $\mathbb K$-valued equations.

The \textbf{$\mathbb K$-valued equation} for the component $u(x,y,s)$ is
$$
-\frac{x}{1-y}\,u_{yy}(x,y,0)
$$
$$
=
\frac{x^2}{(1-y)^2}\,u_{xx}(x,y,0)
-\frac{x(1-y+x)}{(1-y)^2}\,u_{xy}(x,y,0)
+\frac{1-y-x}{(1-y)^2}\,u_y(x,y,0).
$$

The \textbf{$\mathbb K$-valued equation} for the component $v(x,y,s)$ is
\begin{equation}\label{rivn-pryklad-15-v-dali}
v_{yy}(x,y,s)
=
-\frac{x}{1-y}\,v_{xx}(x,y,0)
+
\left(1+\frac{x}{1-y}\right)v_{xy}(x,y,0)
+
\frac{x+y-1}{(1-y)^2}\,v_x(x,y,0).
\end{equation}

\vskip2mm
\textbf{Example 16.} Let $X=\mathcal{O}(\mathbb{C})^2$, where $\mathcal{O}(\mathbb{C})$ is the space of all entire functions of the complex variable $s$,
and let $N$ be a nilpotent operator of index two:
$$
N=
\left(
\begin{array}{cc}
0&1\\
0&0
\end{array}
\right),
\qquad N^2=0.
$$

Elements of the space $X$ have the form
$$
f(s)=
\left(
\begin{array}{c}
f_1(s)\\
f_2(s)
\end{array}
\right),
\qquad
f_1,f_2\in\mathcal{O}(\mathbb{C}).
$$

Set
$$
A(x,y)=xN,\qquad B(x,y)=yN.
$$
Then
$$
A_y=0,\qquad B_x=0,
$$
and hence
$$
A_y=AB_x=0.
$$
Moreover,
$$
B_y=N,\qquad A_x=N,\qquad BB_x=0,
$$
so that
$$
B_y=A_x+BB_x.
$$
Thus, the corresponding algebra is compatible.

Since
$$
A^2=0,\qquad AB=0,\qquad BA_x=0,\qquad AB_x=0,
$$
the \textbf{operator-differential equation} for the component $U(x,y)$ takes the form
\begin{equation}\label{pryklad-16-rivn-U}
xNU_{yy}=0.
\end{equation}

The \textbf{operator-differential equation} for the component $V(x,y)$ has the form
\begin{equation}\label{pryklad-16-rivn-V}
V_{yy}=xNV_{xx}+yNV_{xy}+2NV_x.
\end{equation}

We now pass to the corresponding $\mathbb K$-valued equations.
Apply both sides of equations (\ref{pryklad-16-rivn-U}) and (\ref{pryklad-16-rivn-V}) to an arbitrary vector $f\in X$ and set
$$
u(x,y)=U(x,y)f,\qquad v(x,y)=V(x,y)f.
$$
Let
$$
u(x,y)=
\left(
\begin{array}{c}
u_1(x,y)\\
u_2(x,y)
\end{array}
\right),
\qquad
v(x,y)=
\left(
\begin{array}{c}
v_1(x,y)\\
v_2(x,y)
\end{array}
\right).
$$

Since
$$
N
\left(
\begin{array}{c}
w_1\\
w_2
\end{array}
\right)
=
\left(
\begin{array}{c}
w_2\\
0
\end{array}
\right),
$$
we have
$$
Nu_{yy}
=
\left(
\begin{array}{c}
(u_2)_{yy}\\
0
\end{array}
\right).
$$
Hence, equation (\ref{pryklad-16-rivn-U}) reduces to the system of \textbf{$\mathbb K$-valued equations}
$$
\left\{
\begin{array}{rcl}
x(u_2)_{yy}&=&0,\\
0&=&0.
\end{array}
\right.
$$
Thus, for $x\neq0$, we have
$$
(u_2)_{yy}=0,
$$
whereas the function $u_1$ is not constrained by this equation.

For the function $v$, we have
$$
Nv_{xx}
=
\left(
\begin{array}{c}
(v_2)_{xx}\\
0
\end{array}
\right),
\qquad
Nv_{xy}
=
\left(
\begin{array}{c}
(v_2)_{xy}\\
0
\end{array}
\right),
\qquad
Nv_x
=
\left(
\begin{array}{c}
(v_2)_x\\
0
\end{array}
\right).
$$
Therefore, equation (\ref{pryklad-16-rivn-V}) reduces to the system of \textbf{$\mathbb K$-valued} linear
differential equations with variable coefficients
$$
\left\{
\begin{array}{rcl}
(v_1)_{yy}
&=&
x(v_2)_{xx}+y(v_2)_{xy}+2(v_2)_x,\\[1mm]
(v_2)_{yy}
&=&
0.
\end{array}
\right.
$$

\vskip2mm
\textbf{Example 17.} Let $X=C^\infty[0,\infty)$, and let the Volterra operator $\mathcal{V}$ be defined by
$$
(\mathcal{V}f)(s)=\int_0^s f(\tau)\,d\tau.
$$
Denote
$$
R_y=(I-y\mathcal{V})^{-1}.
$$
Set
$$
P(x,y)=x\mathcal{V} R_y,
\qquad
Q(x,y)=\mathcal{V}.
$$
Since all these operators are functions of the same operator $\mathcal{V}$, they commute. Moreover,
$$
P_y=PP_x,\qquad Q_y=QQ_x=0.
$$
Therefore,
$$
A=-PQ=-x\mathcal V^2R_y,\qquad B=P+Q=x\mathcal V R_y+\mathcal V.
$$
Hence, $A$ and $B$ define a compatible algebra.

The \textbf{operator-differential equation} for the component $U(x,y)$ has the form
$$
-x\mathcal V^2R_yU_{yy}
=
x^2\mathcal V^4R_y^2U_{xx}
-
\left(x\mathcal V^3R_y+x^2\mathcal V^3R_y^2\right)U_{xy}
+
\left(\mathcal V^3R_y-x\mathcal V^3R_y^2\right)U_y.
$$

The \textbf{operator-differential equation} for the component $V(x,y)$ has the form
\begin{equation}\label{pryklad-17-rivn-V}
V_{yy}
=
-x\mathcal V^2R_yV_{xx}
+
\left(x\mathcal VR_y+\mathcal V\right)V_{xy}
+
\left(-\mathcal V^2R_y+x\mathcal V^2R_y^2\right)V_x.
\end{equation}

Here
$$
R_y=(I-y\mathcal V)^{-1}.
$$
In particular,
$$
(R_yf)(s)
=
f(s)+y\int_0^s e^{y(s-\tau)}f(\tau)\,d\tau.
$$

Since the operator-differential equations are rather complicated, we write the corresponding $\mathbb{K}$-valued equation only for the function $v(x,y,s)$.
Apply both sides of the operator-differential equation (\ref{pryklad-17-rivn-V}) to an arbitrary
function $f\in X$ and define
$$
v(x,y,s)=(V(x,y)f)(s).
$$

Thus, the operator-differential equation (\ref{pryklad-17-rivn-V}) for the component $V$ becomes
the following $\mathbb K$-valued integro-differential equation:
$$
v_{yy}(x,y,s)
=
-x\int_0^s
\frac{e^{y(s-\tau)}-1}{y}\,
v_{xx}(x,y,\tau)\,d\tau
+
\int_0^s
\left(1+xe^{y(s-\tau)}\right)
v_{xy}(x,y,\tau)\,d\tau
$$
$$
\quad
-
\int_0^s
\frac{e^{y(s-\tau)}-1}{y}\,
v_x(x,y,\tau)\,d\tau
+
x\int_0^s
(s-\tau)e^{y(s-\tau)}
v_x(x,y,\tau)\,d\tau.
$$

Thus, realizing the structural operators $A$ and $B$ in terms of the
Volterra operator leads to linear
integro-differential equations with variable coefficients for the
$\mathbb K$-valued functions $u(x,y,s)$ and $v(x,y,s)$.

\vskip2mm
\textbf{Example 18.} Let $X=C^\infty(\mathbb{R})$, and let the forward difference operator
$\Delta_h$ be defined by
$$
(\Delta_hf)(s)=f(s+h)-f(s),
\qquad h>0.
$$
Set
$$
A(x,y)=\frac{\lambda}{y+c}\Delta_h,
\qquad
B(x,y)=\frac{a-x}{y+c}I,
\qquad y\neq-c,
$$
where $a,c,\lambda\in\mathbb K$ are constants.

We have
$$
A_x=0,\qquad
A_y=-\frac{\lambda}{(y+c)^2}\Delta_h,
$$
and also
$$
B_x=-\frac1{y+c}I,\qquad
B_y=-\frac{a-x}{(y+c)^2}I.
$$
Therefore,
$$
AB_x
=
-\frac{\lambda}{(y+c)^2}\Delta_h
=
A_y,
$$
and
$$
A_x+BB_x
=
-\frac{a-x}{(y+c)^2}I
=
B_y.
$$
Hence, the operators $A$ and $B$ define a compatible algebra.

The operator-differential equation for the component $V(x,y)$ has the form
$$
V_{yy}
=
\frac{\lambda}{y+c}\Delta_hV_{xx}
+
\frac{a-x}{y+c}V_{xy}
+
\frac{x-a}{(y+c)^2}V_x.
$$

The operator-differential equation for the component $U(x,y)$ has the form
$$
\frac{\lambda}{y+c}\Delta_hU_{yy}
=
\frac{\lambda^2}{(y+c)^2}\Delta_h^2U_{xx}
+
\frac{\lambda(a-x)}{(y+c)^2}\Delta_hU_{xy}
-
\frac{2\lambda}{(y+c)^2}\Delta_hU_y.
$$

We now pass to the corresponding $\mathbb{K}$-valued equations.
Apply both sides of the above equations to an arbitrary function
$f\in X$ and define
$$
u(x,y,s)=(U(x,y)f)(s),\qquad
v(x,y,s)=(V(x,y)f)(s).
$$

Since
$$
(\Delta_hg)(s)=g(s+h)-g(s),
$$
for the function $v$ we obtain the difference-differential equation
$$
v_{yy}(x,y,s)
=
\frac{\lambda}{y+c}
\Biggr(
v_{xx}(x,y,s+h)-v_{xx}(x,y,s)
\Biggr)
$$
$$
\qquad
+
\frac{a-x}{y+c}\,v_{xy}(x,y,s)
+
\frac{x-a}{(y+c)^2}\,v_x(x,y,s).
$$

For the function $u$, we have
$$
\frac{\lambda}{y+c}
\Biggr(
u_{yy}(x,y,s+h)-u_{yy}(x,y,s)
\Biggr)
$$
$$
=
\frac{\lambda^2}{(y+c)^2}
\Biggr(
u_{xx}(x,y,s+2h)
-2u_{xx}(x,y,s+h)
+u_{xx}(x,y,s)
\Biggr)
$$
$$
\quad
+
\frac{\lambda(a-x)}{(y+c)^2}
\Biggr(
u_{xy}(x,y,s+h)-u_{xy}(x,y,s)
\Biggr)
$$
$$
\quad
-
\frac{2\lambda}{(y+c)^2}
\Biggr(
u_y(x,y,s+h)-u_y(x,y,s)
\Biggr).
$$

Thus, the resulting equations are linear
difference-differential equations with variable coefficients.

\vskip2mm
\textbf{Example 19. Inverse problem.} Let $X=C^\infty(\mathbb{R})$ be the space of all infinitely
differentiable functions $f=f(s)$, and let $\alpha=\alpha(s)$ be a fixed smooth function.

Consider the following linear differential equation with variable
coefficients for the $\mathbb{K}$-valued function $v=v(x,y,s)$:
$$
v_{yy}
=
-\frac{\alpha(s)x}{1-\alpha(s)y}\,v_{xx}
+
\left(
1+\frac{\alpha(s)x}{1-\alpha(s)y}
\right)v_{xy}
+
\frac{\alpha(s)\bigl(\alpha(s)x+\alpha(s)y-1\bigr)}
{\bigl(1-\alpha(s)y\bigr)^2}\,v_x.
$$

We determine operator structural coefficients $A(x,y)$ and $B(x,y)$
for which this equation arises as a realization of the
operator-differential equation
$$
V_{yy}
=
AV_{xx}
+
BV_{xy}
+
(2A_x+BB_x)V_x.
$$

Comparing the coefficients of $v_{xx}$ and $v_{xy}$, it is natural to set
$$
A(x,y)
=
M_{-\frac{\alpha(s)x}{1-\alpha(s)y}},
\qquad
B(x,y)
=
M_{1+\frac{\alpha(s)x}{1-\alpha(s)y}}.
$$

Let us verify the coefficient of $v_x$. We have
$$
A_x
=
M_{-\frac{\alpha(s)}{1-\alpha(s)y}},
\qquad
B_x
=
M_{\frac{\alpha(s)}{1-\alpha(s)y}}.
$$
Since the product of multiplication operators is the multiplication operator
by the product of the corresponding functions, we obtain
$$
2A_x+BB_x
=
M_{
-\frac{2\alpha(s)}{1-\alpha(s)y}
+
\left(
1+\frac{\alpha(s)x}{1-\alpha(s)y}
\right)
\frac{\alpha(s)}{1-\alpha(s)y}
}.
$$
After reducing to a common denominator, we obtain
$$
2A_x+BB_x
=
M_{
\frac{\alpha(s)\bigl(\alpha(s)x+\alpha(s)y-1\bigr)}
{\bigl(1-\alpha(s)y\bigr)^2}
}.
$$

Hence, the given $\mathbb{K}$-valued equation is generated by the
operator-differential equation
$$
V_{yy}
=
M_{-\frac{\alpha(s)x}{1-\alpha(s)y}}V_{xx}
+
M_{1+\frac{\alpha(s)x}{1-\alpha(s)y}}V_{xy}
+
M_{
\frac{\alpha(s)\bigl(\alpha(s)x+\alpha(s)y-1\bigr)}
{\bigl(1-\alpha(s)y\bigr)^2}
}V_x.
$$

It remains to verify the compatibility of the structural operators obtained above.
We have
$$
A_y
=
M_{-\frac{\alpha^2(s)x}{\bigl(1-\alpha(s)y\bigr)^2}},
$$
and
$$
AB_x
=
M_{-\frac{\alpha(s)x}{1-\alpha(s)y}}
M_{\frac{\alpha(s)}{1-\alpha(s)y}}
=
M_{-\frac{\alpha^2(s)x}{\bigl(1-\alpha(s)y\bigr)^2}}.
$$
Therefore,
$$
A_y=AB_x.
$$

Moreover,
$$
B_y
=
M_{\frac{\alpha^2(s)x}{\bigl(1-\alpha(s)y\bigr)^2}},
$$
whereas
$$
A_x+BB_x
=
M_{-\frac{\alpha(s)}{1-\alpha(s)y}}
+
M_{
\left(
1+\frac{\alpha(s)x}{1-\alpha(s)y}
\right)
\frac{\alpha(s)}{1-\alpha(s)y}
}
=
M_{\frac{\alpha^2(s)x}{\bigl(1-\alpha(s)y\bigr)^2}}.
$$
Hence,
$$
B_y=A_x+BB_x.
$$

Thus, the operators $A(x,y)$ and $B(x,y)$ obtained above define
a compatible algebra, and the initial $\mathbb{K}$-valued equation is
a realization of the corresponding operator-differential equation
for the component $V(x,y)$.

\vskip2mm
\textbf{Conclusion.}
These examples clearly demonstrate that our theory is not tied to any particular type of equation.
After suitable realizations of the operators $A$ and $B$, the operator-differential equations
(\ref{rivn-V-zmin-operatory}) and (\ref{rivn-U-zmin-operatory})
generate various equations (or systems of equations) for $\mathbb{K}$-valued functions:
\begin{itemize}
  \item linear partial differential equations with variable coefficients;
  \item systems of linear partial differential equations;
  \item difference-differential equations;
  \item integro-differential equations, etc.;
\end{itemize}
and the type of the resulting equation (or system) for a $\mathbb{K}$-valued function is determined by the type of the operators $A$ and $B$ belonging to a commutative
associative algebra of linear operators.

\section{Construction of Solutions to Equations (Systems of Equations) for $\mathbb{K}$-Valued Functions}

In Section 3
we established a relationship between monogenic functions and solutions of the operator-differential equations
(\ref{rivn-V-zmin-operatory}) and (\ref{rivn-U-zmin-operatory}). Namely, all solutions of equations (\ref{rivn-V-zmin-operatory}) and (\ref{rivn-U-zmin-operatory}) can be obtained as components
of monogenic functions.
Then, in Section 4
we demonstrated how different realizations of the operators $A$ and $B$ lead to different types of equations or systems of equations. A natural question now arises: given operator-valued solutions of equations (\ref{rivn-V-zmin-operatory}) and (\ref{rivn-U-zmin-operatory}), how can one construct solutions of the corresponding generated $\mathbb{K}$-valued equations (or systems of equations)? In this section, we answer these questions.

\vskip2mm
\textbf{Theorem 17.}\label{Teor-17-zmin-oper}
Let $X$ be a linear space over the field $\mathbb K$,
let $\mathcal{A}(X)$ be a commutative associative algebra of linear operators,
and let the operator-valued functions $U,V:D\to\mathcal{A}(X)$
possess the required partial derivatives with respect to $x,y$. Suppose that $U$ and $V$ satisfy the operator-differential equations
(\ref{rivn-U-zmin-operatory}) and (\ref{rivn-V-zmin-operatory}), respectively.
Then, for every element $\varphi\in X$, the functions
$$
u(x,y,\xi):=\Biggr(U(x,y)\varphi\Biggr)(\xi),
\qquad
v(x,y,\xi):=\Biggr(V(x,y)\varphi\Biggr)(\xi),
$$
where
$$
\xi=(s,t,\ldots),
$$
satisfy, respectively, the equations
$$
\bigl(Au_{yy}(x,y,\cdot)\bigr)(\xi)
=
\bigl(A^2u_{xx}(x,y,\cdot)\bigr)(\xi)
$$
$$
\qquad
+
\bigl(ABu_{xy}(x,y,\cdot)\bigr)(\xi)
+
\Biggr((2AB_x-BA_x)
u_y(x,y,\cdot)\Biggr)(\xi)
$$
and
$$
v_{yy}(x,y,\xi)
=
\bigl(Av_{xx}(x,y,\cdot)\bigr)(\xi)
+
\bigl(Bv_{xy}(x,y,\cdot)\bigr)(\xi)
$$
\begin{equation}\label{riv-V-teor-17}
\qquad
+
\Biggr((2A_x+BB_x)
v_x(x,y,\cdot)\Biggr)(\xi).
\end{equation}

\textbf{Proof.}
Let $\varphi\in X$ be an arbitrary fixed element and set
$$
u(x,y,\xi):=(U(x,y)\varphi)(\xi),
\qquad
v(x,y,\xi):=(V(x,y)\varphi)(\xi),
$$
where $\xi=(s,t,\ldots)$ denotes the collection of internal variables of the space $X$.

Since the operator-valued functions $U$ and $V$ possess the required partial derivatives with respect to $x,y$, and the element $\varphi$ does not depend on $x,y$, we have
$$
u_x(x,y,\xi)=(U_x(x,y)\varphi)(\xi),
\qquad
u_y(x,y,\xi)=(U_y(x,y)\varphi)(\xi),
$$
and, respectively,
$$
u_{xx}(x,y,\xi)=(U_{xx}(x,y)\varphi)(\xi),
$$
$$
u_{xy}(x,y,\xi)=(U_{xy}(x,y)\varphi)(\xi),
\qquad
u_{yy}(x,y,\xi)=(U_{yy}(x,y)\varphi)(\xi).
$$
Similarly,
$$
v_x(x,y,\xi)=(V_x(x,y)\varphi)(\xi),
$$
$$
v_{xx}(x,y,\xi)=(V_{xx}(x,y)\varphi)(\xi),
\qquad
v_{xy}(x,y,\xi)=(V_{xy}(x,y)\varphi)(\xi),
$$
$$
v_{yy}(x,y,\xi)=(V_{yy}(x,y)\varphi)(\xi).
$$

By assumption, the operator-valued function $U$ satisfies the equation
$$
AU_{yy}
=
A^2U_{xx}
+
ABU_{xy}
+
(2AB_x-BA_x)U_y.
$$
Applying both sides of this equation to the element $\varphi\in X$, we obtain
$$
AU_{yy}\varphi
=
A^2U_{xx}\varphi
+
ABU_{xy}\varphi
+
(2AB_x-BA_x)U_y\varphi.
$$
Evaluating both sides at the point $\xi$ of the internal variable space, we obtain
$$
\bigl(Au_{yy}(x,y,\cdot)\bigr)(\xi)
=
\bigl(A^2u_{xx}(x,y,\cdot)\bigr)(\xi)
$$
$$
\qquad
+
\bigl(ABu_{xy}(x,y,\cdot)\bigr)(\xi)
+
\bigl((2AB_x-BA_x)u_y(x,y,\cdot)\bigr)(\xi).
$$
Hence, the function $u$ satisfies the corresponding $\mathbb K$-valued equation.

Similarly, one proves that the function $v$ satisfies equation (\ref{riv-V-teor-17}).
The theorem is proved.

The converse of Theorem 17 also holds.
\vskip2mm

\textbf{Theorem 18.}
Let $X$ be a linear space over the field $\mathbb{K}$, and let
$A,B:D\to\mathcal{A}(X)$ and
$V:D\to\mathcal{A}(X)$ possess all the required partial derivatives.
Suppose that, for every element $\varphi\in X$, the function
$v_\varphi(x,y):=V(x,y)\varphi$ satisfies the equation
$$
(v_\varphi)_{yy}
=
A(v_\varphi)_{xx}
+B(v_\varphi)_{xy}
+\bigl(A_x+B_y+BA_y\bigr)(v_\varphi)_x.
$$
Then the operator-valued function $V$ satisfies the operator-differential equation
$$
V_{yy}
=
A V_{xx}
+B V_{xy}
+\bigl(A_x+B_y+BA_y\bigr)V_x.
$$

\textbf{Proof.}
For an arbitrary $\varphi\in X$, we have
$$
(v_\varphi)_x=V_x\varphi,\qquad
(v_\varphi)_{xx}=V_{xx}\varphi,\qquad
(v_\varphi)_{xy}=V_{xy}\varphi,\qquad
(v_\varphi)_{yy}=V_{yy}\varphi.
$$
Therefore,
$$
\Biggr(
V_{yy}
-A V_{xx}
-B V_{xy}
-\bigl(A_x+B_y+BA_y\bigr)V_x
\Biggr)\varphi=0
$$
for every $\varphi\in X$.
Hence,
$$
V_{yy}
-A V_{xx}
-B V_{xy}
-\bigl(A_x+B_y+BA_y\bigr)V_x=0,
$$
which proves the assertion.

An analogous theorem holds for the function
$u_\varphi(x,y):=U(x,y)\varphi$.

Let us consider some examples.

\vskip2mm
\textbf{Example 20.}
Apply Theorem 17
to Example 13. In this case,
$$
A(x,y)=-\frac{x}{1-y}I,\qquad
B(x,y)=\left(1+\frac{x}{1-y}\right)I,\qquad y\ne1.
$$

The corresponding operator-differential equation for the component $U$ has the form
\begin{equation}\label{rivn-pryklad-13-U-dali}
-\frac{x}{1-y}\,U_{yy}
=
\frac{x^2}{(1-y)^2}\,U_{xx}
-\frac{x(1-y+x)}{(1-y)^2}\,U_{xy}
+\frac{1-y-x}{(1-y)^2}\,U_y,
\end{equation}
whereas the equation for the component $V$ has the form
\begin{equation}\label{rivn-pryklad-13-V-dali}
V_{yy}
=
-\frac{x}{1-y}\,V_{xx}
+\left(1+\frac{x}{1-y}\right)V_{xy}
+\frac{x+y-1}{(1-y)^2}\,V_x.
\end{equation}

After realizing the identity operator $I$ on the space $X=C^\infty(\mathbb{R})$, we obtain the corresponding $\mathbb{K}$-valued equations
\begin{equation}\label{rivn-pryklad-13-u-dali}
-\frac{x}{1-y}\,u_{yy}
=
\frac{x^2}{(1-y)^2}\,u_{xx}
-\frac{x(1-y+x)}{(1-y)^2}\,u_{xy}
+\frac{1-y-x}{(1-y)^2}\,u_y
\end{equation}
and
\begin{equation}\label{rivn-pryklad-13-v-dali}
v_{yy}
=
-\frac{x}{1-y}\,v_{xx}
+\left(1+\frac{x}{1-y}\right)v_{xy}
+\frac{x+y-1}{(1-y)^2}\,v_x.
\end{equation}

Consider the operator-valued monogenic function
$$
\Phi(Z)=Z^{\star2}=UE_1+VE_2,
\qquad Z=xE_1+yE_2.
$$
Since
$$
E_2\star E_2=A(x,y)E_1+B(x,y)E_2,
$$
we have
\begin{equation}\label{rivn-U-V-pryklady}
U=x^2I+y^2A,\qquad
V=2xyI+y^2B.
\end{equation}
Substituting the operator-valued functions $A$ and $B$, we obtain
$$
U(x,y)=
\left(x^2-\frac{xy^2}{1-y}\right)I
$$
and
$$
V(x,y)=
\left(
2xy+y^2+\frac{xy^2}{1-y}
\right)I.
$$

Let $\varphi\in C^\infty(\mathbb{R})$. Then, for the component $U$, we have
$$
u(x,y,s):=\bigl(U(x,y)\varphi\bigr)(s)
=
\left(x^2-\frac{xy^2}{1-y}\right)(I\varphi)(s),
$$
that is,
\begin{equation}\label{rozviazok-pryklad-13-u}
u(x,y,s)
=
\left(x^2-\frac{xy^2}{1-y}\right)\varphi(s).
\end{equation}

For the component $V$, we obtain
$$
v(x,y,s):=\bigl(V(x,y)\varphi\bigr)(s),
$$
that is,
\begin{equation}\label{rozviazok-pryklad-13-v}
v(x,y,s)=
\left(
2xy+y^2+\frac{xy^2}{1-y}
\right)\varphi(s).
\end{equation}

Hence, by Theorem 17,
the function (\ref{rozviazok-pryklad-13-u}), where $\varphi$ is an arbitrary function from $C^\infty(\mathbb{R})$, is a $\mathbb{K}$-valued solution of equation (\ref{rivn-pryklad-13-u-dali}), whereas the function (\ref{rozviazok-pryklad-13-v}) is a $\mathbb{K}$-valued solution of equation (\ref{rivn-pryklad-13-v-dali}).

Indeed, since the function $\varphi(s)$ does not depend on $x$ and $y$, after substituting the functions (\ref{rozviazok-pryklad-13-u}) and (\ref{rozviazok-pryklad-13-v}) into the corresponding equations, the factor $\varphi(s)$ can be taken outside all partial derivatives with respect to $x$ and $y$, and the identities can be verified by direct calculation.

Similarly, every operator-valued monogenic function in the algebra under consideration, after acting on an arbitrary element $\varphi\in C^\infty(\mathbb{R})$, generates explicit $\mathbb{K}$-valued solutions of the corresponding equations (\ref{rivn-pryklad-13-u-dali}) and (\ref{rivn-pryklad-13-v-dali}).

\vskip2mm
\textbf{Example 21.}
Apply Theorem 17
to Example 14.1.
In this case,
$$
A(x,y)=M_{-\frac{\alpha(s)x}{1-\alpha(s)y}},
\qquad
B(x,y)=M_{1+\frac{\alpha(s)x}{1-\alpha(s)y}},
\qquad
1-\alpha(s)y\ne0.
$$

The corresponding operator-differential equation for the component $V$ has the form
\begin{equation}\label{rivn-pryklad-14-V-dali}
V_{yy}
=
M_{-\frac{\alpha(s)x}{1-\alpha(s)y}}\,V_{xx}
+
M_{1+\frac{\alpha(s)x}{1-\alpha(s)y}}\,V_{xy}
+
M_{\frac{\alpha(s)\bigl(\alpha(s)x+\alpha(s)y-1\bigr)}
{\bigl(1-\alpha(s)y\bigr)^2}}\,V_x.
\end{equation}

After realizing the multiplication operators on the space
$X=C^\infty(\mathbb{R})$, we obtain the corresponding
$\mathbb{K}$-valued differential equation
\begin{equation}\label{rivn-pryklad-14-v-dali}
v_{yy}
=
-\frac{\alpha(s)x}{1-\alpha(s)y}\,v_{xx}
+
\left(
1+\frac{\alpha(s)x}{1-\alpha(s)y}
\right)v_{xy}
+
\frac{\alpha(s)\bigl(\alpha(s)x+\alpha(s)y-1\bigr)}
{\bigl(1-\alpha(s)y\bigr)^2}\,v_x.
\end{equation}

For the operator-valued monogenic function
$\Phi(Z)=Z^{\star2}$, the components have the form (\ref{rivn-U-V-pryklady}).

Substituting the operator-valued function $B$, we obtain
$$
V(x,y)
=
2xyI+
y^2M_{1+\frac{\alpha(s)x}{1-\alpha(s)y}}.
$$
Since $2xyI=M_{2xy}$, we have
$$
V(x,y)
=
M_{2xy+y^2+\frac{\alpha(s)xy^2}{1-\alpha(s)y}}.
$$

Let $\varphi\in C^\infty(\mathbb{R})$. Then
$$
v(x,y,s):=\bigl(V(x,y)\varphi\bigr)(s),
$$
that is,
\begin{equation}\label{rozviazok-pryklad-14-v}
v(x,y,s)
=
\left(
2xy+y^2+
\frac{\alpha(s)xy^2}{1-\alpha(s)y}
\right)\varphi(s).
\end{equation}

Hence, by Theorem 17,
the function
(\ref{rozviazok-pryklad-14-v}), where $\varphi$ is an arbitrary function from
$C^\infty(\mathbb{R})$, is a $\mathbb{K}$-valued solution of equation
(\ref{rivn-pryklad-14-v-dali}). This can easily be verified directly.

Similarly, every operator-valued monogenic function in the algebra under consideration, after acting on an arbitrary element
$\varphi\in C^\infty(\mathbb{R})$, generates explicit $\mathbb{K}$-valued solutions of equation (\ref{rivn-pryklad-14-v-dali}).

\vskip2mm
\textbf{Example 22.}
Apply Theorem 17
to Example 15. In this case,
$$
A(x,y)=-\frac{x}{1-y}P,\qquad
B(x,y)=\left(1+\frac{x}{1-y}\right)P,
\qquad y\ne1.
$$

The corresponding operator-differential equation for the component $V$ has the form
(\ref{rivn-pryklad-15-V-dali}). After realizing the projection operator $P$ on the space
$X=C^\infty(\mathbb{R})$, where $(Pf)(s)=f(0)$,
we obtain the corresponding $\mathbb{K}$-valued equation (\ref{rivn-pryklad-15-v-dali}):
$$
v_{yy}(x,y,s)
=
-\frac{x}{1-y}\,v_{xx}(x,y,0)
+
\left(1+\frac{x}{1-y}\right)v_{xy}(x,y,0)
+
\frac{x+y-1}{(1-y)^2}\,v_x(x,y,0).
$$

Consider the operator-valued monogenic function
$$
\Phi(Z)=Z^{\star2}=UE_1+VE_2,
$$
where
$V=2xyI+y^2B$.

Substituting the operator-valued function $B$, we obtain
$$
V(x,y)
=
2xyI+
y^2\left(1+\frac{x}{1-y}\right)P.
$$

Let $\varphi\in C^\infty(\mathbb{R})$. Then
$$
v(x,y,s):=\bigl(V(x,y)\varphi\bigr)(s).
$$
Since
$$
(I\varphi)(s)=\varphi(s),
\qquad
(P\varphi)(s)=\varphi(0),
$$
we obtain
\begin{equation}\label{rozviazok-pryklad-15-v}
v(x,y,s)
=
2xy\,\varphi(s)
+
y^2\left(1+\frac{x}{1-y}\right)\varphi(0).
\end{equation}

Hence, by Theorem 17,
the function (\ref{rozviazok-pryklad-15-v}), where
$\varphi$ is an arbitrary function from
$C^\infty(\mathbb{R})$, is a
$\mathbb{K}$-valued solution of equation
(\ref{rivn-pryklad-15-v-dali}).

\vskip2mm
\textbf{Example 23.}
We now apply Theorem 17
to Example 16. In this case,
$$
A(x,y)=xN,\qquad B(x,y)=yN,
$$
where
$$
N=
\left(
\begin{array}{cc}
0&1\\
0&0
\end{array}
\right),
\qquad N^2=0.
$$

The corresponding operator-differential equation for the component $V$ has the form
\begin{equation}\label{pryklad-16-rivn-V-dali}
V_{yy}=xNV_{xx}+yNV_{xy}+2NV_x.
\end{equation}

After applying both sides of equation (\ref{pryklad-16-rivn-V-dali}) to an arbitrary vector
$\varphi\in\mathbb{K}^2$ and setting
$$
v(x,y):=V(x,y)\varphi=
\left(
\begin{array}{c}
v_1(x,y)\\
v_2(x,y)
\end{array}
\right),
$$
we obtain the system of $\mathbb{K}$-valued equations
\begin{equation}\label{pryklad-16-systema-v-dali}
\left\{
\begin{array}{rcl}
(v_1)_{yy}&=&x(v_2)_{xx}+y(v_2)_{xy}+2(v_2)_x,\\
(v_2)_{yy}&=&0.
\end{array}
\right.
\end{equation}

Consider the operator-valued monogenic function
$$
\Phi(Z)=Z^{\star3}=UE_1+VE_2,
\qquad
Z=xE_1+yE_2.
$$
Here
$$
V(x,y)=
3x^2yI+4xy^3N.
$$

Let
$$
\varphi(s)=
\left(
\begin{array}{c}
\varphi_1(s)\\
\varphi_2(s)
\end{array}
\right)
\in\mathcal{O}(\mathbb{C})^2.
$$
Then
$$
v(x,y,s):=V(x,y)\varphi(s).
$$
Since
$$
N\varphi(s)=
\left(
\begin{array}{c}
\varphi_2(s)\\
0
\end{array}
\right),
$$
we have
$$
v(x,y,s)=
3x^2y
\left(
\begin{array}{c}
\varphi_1(s)\\
\varphi_2(s)
\end{array}
\right)
+
4xy^3
\left(
\begin{array}{c}
\varphi_2(s)\\
0
\end{array}
\right).
$$
Hence,
\begin{equation}\label{rozviazok-pryklad-16-v}
v(x,y,s)=
\left(
\begin{array}{c}
3x^2y\,\varphi_1(s)+4xy^3\,\varphi_2(s)\\[1mm]
3x^2y\,\varphi_2(s)
\end{array}
\right).
\end{equation}

Thus,
\begin{equation}\label{rozviazok-pryklad-16-v1}
v_1(x,y,s)=3x^2y\,\varphi_1(s)+4xy^3\,\varphi_2(s),
\end{equation}
and
\begin{equation}\label{rozviazok-pryklad-16-v2}
v_2(x,y,s)=3x^2y\,\varphi_2(s).
\end{equation}

Hence, by Theorem 17, the function (\ref{rozviazok-pryklad-16-v}), where
$\varphi_1,\varphi_2\in\mathcal{O}(\mathbb{C})$ are arbitrary entire functions, is a
$\mathbb{K}^2$-valued solution of equation (\ref{pryklad-16-rivn-V-dali}), while its components
(\ref{rozviazok-pryklad-16-v1}) and (\ref{rozviazok-pryklad-16-v2}) are $\mathbb{K}$-valued solutions of the system
(\ref{pryklad-16-systema-v-dali}).

Similarly, every operator-valued monogenic function in the algebra under consideration, after acting on an arbitrary element
$\varphi\in\mathcal O(\mathbb C)^2$, generates explicit $\mathbb{K}^2$-valued solutions of equation
(\ref{pryklad-16-rivn-V-dali}) and the corresponding $\mathbb{K}$-valued solutions of system (\ref{pryklad-16-systema-v-dali}).

\vskip2mm
\textbf{Example 24.}
We now apply Theorem 17
to Example 17. In this case,
$$
A(x,y)=-x\mathcal{V}^2R_y,
\qquad
B(x,y)=x\mathcal{V}R_y+\mathcal{V},
\qquad
R_y=(I-y\mathcal{V})^{-1}.
$$

The corresponding operator-differential equation for the component $V$ has the form
\begin{equation}\label{pryklad-24-rivn-V}
V_{yy}
=
-x\mathcal{V}^2R_yV_{xx}
+
\left(x\mathcal{V}R_y+\mathcal{V}\right)V_{xy}
+
\left(-\mathcal{V}^2R_y+x\mathcal{V}^2R_y^2\right)V_x.
\end{equation}

After realizing the operators on the space
$X=C^\infty[0,\infty)$, we obtain the corresponding $\mathbb{K}$-valued integro-differential equation
$$
v_{yy}(x,y,s)
=
-x\int_0^s
\frac{e^{y(s-\tau)}-1}{y}\,
v_{xx}(x,y,\tau)d\tau
+
\int_0^s \left(1+xe^{y(s-\tau)}\right)
v_{xy}(x,y,\tau)d\tau
$$
\begin{equation}\label{pryklad-24-rivn-v}
-
\int_0^s
\frac{e^{y(s-\tau)}-1}{y}\,
v_x(x,y,\tau)d\tau
+
x\int_0^s
(s-\tau)e^{y(s-\tau)}
v_x(x,y,\tau)d\tau .
\end{equation}

For the operator-valued monogenic function
$\Phi(Z)=Z^{\star2}$, the components have the form (\ref{rivn-U-V-pryklady}).

Substituting the operator-valued function $B$, we obtain
$$
V(x,y)
=
2xyI+y^2\left(x\mathcal{V}R_y+\mathcal{V}\right),
$$
that is,
\begin{equation}\label{pryklad-24-operator-V}
V(x,y)
=
2xyI+xy^2\mathcal{V}R_y+y^2\mathcal{V}.
\end{equation}

Let $f\in C^\infty[0,\infty)$. Then
$$
v(x,y,s):=\bigl(V(x,y)f\bigr)(s).
$$
Since
$(\mathcal{V}f)(s)=\int_0^s f(\tau)d\tau$
and
$(\mathcal{V}R_yf)(s)=\int_0^s e^{y(s-\tau)}f(\tau)d\tau$,
from (\ref{pryklad-24-operator-V}) we obtain
\begin{equation}\label{pryklad-24-rozviazok-v}
v(x,y,s)
=
2xyf(s)
+
xy^2\int_0^s e^{y(s-\tau)}f(\tau)d\tau
+
y^2\int_0^s f(\tau)d\tau .
\end{equation}

Hence, by Theorem 17,
the function (\ref{pryklad-24-rozviazok-v}), where $f$ is an arbitrary function from
$C^\infty[0,\infty)$, is a $\mathbb{K}$-valued solution of the integro-differential equation
(\ref{pryklad-24-rivn-v}).

Thus, the operator-valued monogenic function
$\Phi(Z)=Z^{\star2}$, after acting on an arbitrary element
$f\in C^\infty[0,\infty)$, generates an explicit infinite-dimensional family of $\mathbb{K}$-valued solutions of equation
(\ref{pryklad-24-rivn-v}).

Similarly, every operator-valued monogenic function in the algebra under consideration, after acting on an arbitrary element
$f\in C^\infty[0,\infty)$,
generates explicit infinite-dimensional families of $\mathbb{K}$-valued solutions of the integro-differential equation (\ref{pryklad-24-rivn-v}).

\vskip2mm
\textbf{Example 25.}
Apply Theorem 17
to Example 18. In this case,
$$
A(x,y)=\frac{\lambda}{y+c}\Delta_h,
\qquad
B(x,y)=\frac{a-x}{y+c}I,
\qquad y\ne-c.
$$

The corresponding operator-differential equation for the component $V$ has the form
\begin{equation}\label{pryklad-25-rivn-V}
V_{yy}
=
\frac{\lambda}{y+c}\Delta_hV_{xx}
+
\frac{a-x}{y+c}V_{xy}
+
\frac{x-a}{(y+c)^2}V_x.
\end{equation}

After realizing the forward difference operator on the space
$X=C^\infty(\mathbb{R})$, we obtain the corresponding $\mathbb{K}$-valued difference-differential equation
$$
v_{yy}(x,y,s)
=
\frac{\lambda}{y+c}
\Biggr(
v_{xx}(x,y,s+h)-v_{xx}(x,y,s)
\Biggr)
$$
\begin{equation}\label{pryklad-25-rivn-v}
+
\frac{a-x}{y+c}\,v_{xy}(x,y,s)
+
\frac{x-a}{(y+c)^2}\,v_x(x,y,s).
\end{equation}

First, consider the operator-valued monogenic function
$\Phi(Z)=Z^{\star2}$. Its component $V$ has the form (\ref{rivn-U-V-pryklady}).
Substituting the operator-valued function $B$, we obtain
$$
V(x,y)
=
2xyI+
\frac{(a-x)y^2}{y+c}I
=
\left(
2xy+\frac{(a-x)y^2}{y+c}
\right)I.
$$

Let $f\in C^\infty(\mathbb{R})$. Then
$$
v_2(x,y,s):=\bigl(V(x,y)f\bigr)(s),
$$
and hence
\begin{equation}\label{pryklad-25-rozviazok-v-2}
v_2(x,y,s)
=
\left(
2xy+\frac{(a-x)y^2}{y+c}
\right)f(s).
\end{equation}

Hence, by Theorem 17,
the function (\ref{pryklad-25-rozviazok-v-2}), where $f$ is an arbitrary function from
$C^\infty(\mathbb{R})$, is a $\mathbb{K}$-valued solution of the difference-differential equation
(\ref{pryklad-25-rivn-v}).

Now consider the operator-valued monogenic function
$\Phi(Z)=Z^{\star3}$.
The component $V$ of the function $Z^{\star3}$ has the form
$$
V(x,y)
=
3x^2yI+y^3A+3xy^2B+y^3B^2.
$$

Substituting the operators $A$ and $B$, we obtain
$$
V(x,y)
=
3x^2yI
+
\frac{\lambda y^3}{y+c}\Delta_h
+
\frac{3xy^2(a-x)}{y+c}I
+
\frac{y^3(a-x)^2}{(y+c)^2}I.
$$
Hence,
\begin{equation}\label{pryklad-25-operator-V-3}
V(x,y)
=
\left(
3x^2y
+
\frac{3xy^2(a-x)}{y+c}
+
\frac{y^3(a-x)^2}{(y+c)^2}
\right)I
+
\frac{\lambda y^3}{y+c}\Delta_h.
\end{equation}

For an arbitrary function $f\in C^\infty(\mathbb{R})$, set
$$
v_3(x,y,s):=\bigl(V(x,y)f\bigr)(s).
$$
Since
$(\Delta_hf)(s)=f(s+h)-f(s)$,
from (\ref{pryklad-25-operator-V-3}) we obtain
$$
v_3(x,y,s)
=
\left(
3x^2y
+
\frac{3xy^2(a-x)}{y+c}
+
\frac{y^3(a-x)^2}{(y+c)^2}
\right)f(s)
$$
\begin{equation}\label{pryklad-25-rozviazok-v-3}
+
\frac{\lambda y^3}{y+c}
\Biggr(f(s+h)-f(s)\Biggr).
\end{equation}

Thus, by Theorem 17,
the function (\ref{pryklad-25-rozviazok-v-3}), where $f$ is an arbitrary function from
$C^\infty(\mathbb{R})$, is also a $\mathbb{K}$-valued solution of the difference-differential equation
(\ref{pryklad-25-rivn-v}).

Similarly, every operator-valued monogenic function in the algebra under consideration, after acting on an arbitrary element
$f\in C^\infty(\mathbb{R})$, generates explicit $\mathbb{K}$-valued solutions of the difference-differential equation
(\ref{pryklad-25-rivn-v}).

\section{Concluding Remarks}

1. In a compatible algebra $\mathfrak{A}_2(\mathcal{A}(X),\star;x,y)$, equation (\ref{rivn-V-zmin-operatory}) can be written in the form
\begin{equation}\label{zakl-zmin-oper-V}
V_{yy}=A V_{xx}+B V_{xy}+(A_x+B_y)V_x.
\end{equation}
Similarly, in a compatible algebra $\mathfrak{A}_2(\mathcal{A}(X),\star;x,y)$, if the operator $A$ is invertible, equation (\ref{rivn-U-zmin-operatory}) can be written in the form
\begin{equation}\label{zakl-zmin-oper-U}
U_{yy}=A U_{xx}+B U_{xy}+\left(2B_x-BA^{-1}A_x\right)U_y.
\end{equation}

2. If we want monogenic functions to be associated, instead of the operator-differential equations (\ref{rivn-V-zmin-operatory}) and (\ref{rivn-U-zmin-operatory}), with an operator-differential equation involving a larger number of independent variables $x$, $y$, $z$,\ldots, then it is necessary to consider a commutative algebra with operator-valued structural constants of higher dimension.

3. If an operator-differential equation has a different form, then instead of the variable $Z=xE_1+yE_2$, one has to consider another variable
$Z=x\widetilde{E}_1+y\widetilde{E}_2$,
where $\widetilde{E}_1,\widetilde{E}_2$ are not necessarily basis elements but certain other elements of the commutative operator algebra.

4. To investigate solutions of the operator-differential equations (\ref{rivn-V-zmin-operatory}) and (\ref{rivn-U-zmin-operatory}), or of other higher-dimensional operator-differential equations (and, consequently, of all corresponding equations and systems for $\mathbb{K}$-valued functions), it becomes necessary to develop an analogue of commutative hypercomplex analysis in the operator algebras introduced above.

\subsection*{Acknowledgment}
 This work was supported by a grant from the Simons Foundation 
(SFI-PD-Ukraine-00014586,V.S.Sh.). The author was supported by budget program “Support for the development of priority areas of research” 
(KPKVK 6541230).

\renewcommand{\refname}{References}

\end{document}